%% file: pderivedseries4_1_07.tex
\documentclass{amsart}
\usepackage{amsmath,amsthm,amscd,amssymb}

\usepackage{color}
\usepackage{amssymb}

\input amsmath-defs.tex
\newcommand{\db}{\bar{\delta}}
\newcommand{\R}{\operatorname{r}}
\newcommand{\Z}{\mathbb{Z}}
\newcommand{\K}{\mathbb{K}}
\newcommand{\rk}{\operatorname{rank}}
\newcommand{\xp}{X^\prime}
\newcommand{\onto}{\twoheadrightarrow}
\newcommand{\Ext}{\operatorname{Ext}}
\newcommand{\red}{\textcolor{red}}

\newtheorem{thm}{Theorem}[section]
\newtheorem{lem}[thm]{Lemma}
\newtheorem{prop}[thm]{Proposition}
\newtheorem{cor}[thm]{Corollary}

\theoremstyle{definition}

\newtheorem{defn}[thm]{Definition}

\newcommand{\bz}{\mathbb{Z}}
\newcommand{\bq}{\mathbb{Q}}

\newcommand{\hz}[1]{\ensuremath{H_#1(-;\mathbb{Z})}}
\newcommand{\und}{\ensuremath{\un{\ \ }}}
\newcommand{\sll}{\ensuremath{S^3\backslash L}}
\newcommand{\sk}{\ensuremath{S^3\backslash K}}

\newcommand{\hq}[1]{\ensuremath{H_#1(-;\mathbb{Q})}}
\newcommand{\der}{\operatorname{Der}}
\newcommand{\id}{\operatorname{id}}
\newcommand{\ho}{\operatorname{Hom}}
\newcommand{\rank}{\operatorname{rank}}
\newcommand{\image}{\operatorname{image}}
\usepackage{pstricks}
\usepackage{pst-node}

\newcommand{\sss}{\scriptscriptstyle}
\newcommand{\2}[1]{\ensuremath{^{\sss (#1)}}}
\newcommand{\np}{\ensuremath{^{\sss (n+1)}}}
\newcommand{\nph}[1]{\ensuremath{#1^{\np}_{\sss H}}}
\newcommand{\nm}{\ensuremath{^{\sss (n-1)}}}
\newcommand{\gn}{\ensuremath{G^{\sss (n)}_{\sss H}}}
\newcommand{\gnp}{\ensuremath{G^{\sss (n+1)}_{\sss H}}}
\newcommand{\gone}{\ensuremath{G^{\sss (1)}_{\sss H}}}
\newcommand{\an}{\ensuremath{A^{\sss (n)}}}
\newcommand{\bn}{\ensuremath{B^{\sss (n)}}}
\newcommand{\cn}{\ensuremath{C^{\sss (n)}}}
\newcommand{\anp}{\ensuremath{A^{\sss (n+1)}}}
\newcommand{\bnp}{\ensuremath{B^{\sss (n+1)}}}
\newcommand{\aln}{\ensuremath{A_{p,n}}}
\newcommand{\bln}{\ensuremath{B_{p,n}}}
\newcommand{\alnp}{\ensuremath{A_{p,n+1}}}
\newcommand{\blnp}{\ensuremath{B_{p,n+1}}}

\title{Homology and Derived p-Series of Groups}
\author{Tim Cochran and Shelly Harvey$^\dag$}

\address{Rice University, Houston, Texas, 77005-1892}
\email{cochran@rice.edu, shelly@rice.edu}

\thanks{$^{\dag}$Both authors were partially supported by the National
Science Foundation. The second author was partially supported
by the Sloan Foundation}

\begin{document}

\begin{abstract} We prove that groups that are $\mathbb{Z}_p$-homology equivalent are isomorphic modulo any term of their $\textbf{derived p-series}$, in precise analogy to Stallings' $1963$ result for the $p-$lower-central series. In fact we prove a stronger theorem that is analogous to Dwyer's extensions of Stallings' results. Similarly, spaces that are $\mathbb{Z}_p$-homology equivalent have isomorphic fundamental groups modulo any term of their $p$-derived series. Various authors have related the ranks of the successive quotients of the $p$-lower central series and of the derived $p$-series of the fundamental group of a $3$-manifold $M$ to the volume of $M$, to whether certain subgroups of $\pi_1(M)$ are free, whether finite index subgroups of $\pi_1(M)$ map onto non-abelian free groups, and to whether finite covers of $M$ are ``large'' in various other senses. 
\end{abstract}

\maketitle

\section{Introduction}\label{intro}

In 1965 J. Stallings showed that if a group homomorphism induces a homology equivalence (with coefficients in $R=\mathbb{Z}$, $\mathbb{Q}$, or $\mathbb{Z}_p)$ then it induces an isomorphism between the groups modulo any term of their respective \textbf{$R$-lower central series} (~\cite{St}, see also ~\cite[Theorem 9.11]{HS}). This was generalized in significant ways by A. Bousfield and W. Dwyer ~\cite{Dw}~\cite{B}. These results have had a significant impact in topology. In the present paper we establish the analogues of Stallings' theorem and Dwyer's theorem for the \textbf{derived p-series}, that is to say, the precise analogue for the derived series in the case $R=\mathbb{Z}_p$. We note that for $R=\mathbb{Z}$ and $R=\mathbb{Q}$,  the precise analogues of Stallings' theorems for the \emph{derived series} are false. However, in 2004, the second author introduced a new series, the \emph{torsion-free derived series}, with which the authors proved analogues of the theorems of Stallings and Dwyer  for $R=\mathbb{Q}$ ~\cite{Ha2}\cite{CH1}\cite{CH2}. By contrast, the case the case $R=\mathbb{Z}_p$ does not exhibit such complications.

The case $R=\mathbb{Z}_p$ is especially important in the study of $3$-manifolds, whose fundamental group controls most of their topology and where questions about the behavior of the homology of finite covers abound. Various authors have related the ranks of the successive quotients of the $p$-lower central series and of the derived $p$-series to volumes of $3$-manifolds ~\cite{CuSh}, whether certain subgroups of $\pi_1$ of a $3$-manifold are free ~\cite{SW}~\cite{CH1}, whether finite index subgroups of $\pi_1(M)$ map onto non-abelian free groups ~\cite{Ho}~\cite{Lu1}\cite{La1}\cite{La2}\cite{La3}, and whether or not finite covers of $M$ are ``large'' in various other senses.

Recall that the $\textbf{lower central series}$ of
$G$, denoted $\{G_{n}\}$, is inductively defined by 
$$
G_1=G ,~~G_{n +1}=[G_{n},G].
$$
Fix a prime $p$. Our convention is that $\mathbb{Z}_p$ denotes the integers modulo $p$. The $\textbf{p-lower central series}$, $\{G_{p,n}\}$ is the fastest descending central series whose successive quotients are $\mathbb{Z}_p$-vector spaces ~\cite{St}. It is defined by
$$
G_{p,1}=G, ~~G_{p,n +1}= (G_{p,n})^p[G_{p,n},G].
$$
Recall Stallings' result:
\begin{thm}[Stallings' $\mathbb{Z}_p$ Theorem {\cite[Theorem~3.4]{St}}]\label{thm:StallZp} Let $\phi:A\to B$ be a
homomorphism that induces an isomorphism on $H_1(-;\mathbb{Z}_p)$ and an
epimorphism on $H_2(-;\mathbb{Z}_p)$. Then, for each $n$, $\phi$
induces an isomorphism $A/A_{p,n}\cong B/B_{p,n}$.
\end{thm}
\noindent Consequently, as was shown by Bousfield

\begin{cor}\label{cor:prop}\cite[Lemma 3.7,Corollary 3.15, Remark 12.2]{B} Under the hypotheses of Theorem~\ref{thm:StallZp}, if dim~$H_1(A;\mathbb{Z}_p)$ is finite, then $\phi$ induces an isomorphism between the pro-$p$-completions $\hat{A}\to\hat{B}$.
\end{cor}
 
The $\textbf{derived p-series}$ (sometimes called the $p$-derived series) denoted $\{G^{\sss (n)}\}$ ($\textbf{warning: the $p$ will be suppressed throughout this paper}$), is defined by
$$
G\20=G,~~ G^{\sss (n +1)}=[G^{\sss (n)},G^{\sss (n)}](G^{\sss (n)})^p.
$$
It is the fastest descending normal series whose successive quotients are $\mathbb{Z}_p$-vector spaces. We also set
$$
G^{(\omega)}=~\cap_{n=1}^\infty G^{(n)}.
$$
Our first result is then:

\newtheorem*{cor:main}{Corollary~\ref{cor:main}}
\begin{cor:main}  Let $A$ and $B$ be finitely-generated groups. If $\phi:A\to B$ induces
an isomorphism (respectively monomorphism) on $H_1(- ;\mathbb{Z}_p)$ and an
epimorphism on $H_2(- ;\mathbb{Z}_p)$, then for each finite $n$,
it induces an isomorphism (respectively monomorphism) $A/\an\to B/\bn$, and a monomorphism $A/A^{(\omega)} \subset B/B^{(\omega)}$. 
\end{cor:main}

Actually we show that this theorem can be derived easily from Stallings' theorem for the $p$-lower central series. Here it is also obtained as a corollary of our main theorem ( Theorem~\ref{thm:dwyermain} below) that is the analogue of Dwyer's extensions of Stallings' theorems.

We have the standard applications.

\newtheorem*{cor:homologyequiv}{Corollary~\ref{cor:homologyequiv}}
\begin{cor:homologyequiv} Suppose $Y$ and $X$ are path-connected $CW$-complexes with $\pi_1(Y)$ and $\pi_1(X)$ finitely generated. If $f:Y\to X$ is a continuous map that induces
an isomorphism (respectively monomorphism) on $H_1(- ;\mathbb{Z}_p)$ and an
epimorphism on $H_2(- ;\mathbb{Z}_p)$, then for each finite $n$,
it induces an isomorphism (respectively monomorphism)
$$
f_*:\pi_1(Y)/\pi_1(Y)^{(n)}\to \pi_1(X)/\pi_1(X)^{(n)}.
$$
\end{cor:homologyequiv}

\newtheorem*{cor:homologycob}{Corollary~\ref{cor:homologycob}}
\begin{cor:homologycob} If $M$ and $N$ are compact $\mathbb{Z}_p$-homology cobordant manifolds then for each $n$
$$
\pi_1(M)/\pi_1(M)^{(n)}\cong \pi_1(N)/\pi_1(N)^{(n)}.
$$
\end{cor:homologycob}

Note that maps $f:N\to M$ between closed orientable $3$-manifolds are particularly nice because if $f_*$ is injective on $H_1(-;\mathbb{Z}_p)$, then it follows from Poincar\'{e} duality that $f_*$ is surjective on $H_2(-;\mathbb{Z}_p)$.

\newtheorem*{cor:hyperbolic}{Corollary~\ref{cor:hyperbolic}}
\begin{cor:hyperbolic} For any closed orientable $3$-manifold $N$,
there exists a hyperbolic $3$-manifold $M$ such that for each $n$
$$
\pi_1(M)/\pi_1(M)^{(n)}\cong \pi_1(N)/\pi_1(N)^{(n)}
$$
and
$$
\pi_1(M)^{(n)}/\pi_1(M)^{(n+1)}\cong \pi_1(N)^{(n)}/\pi_1(N)^{(n+1)}.
$$
In particular the growth rate of the successive quotients of the derived $p$-series of an arbitrary $3$-manifold is the same as that of a hyperbolic $3$-manifold with the same integral homology. Moreover $\pi_1(M)$ and $\pi_1(N)$ have isomorphic pro-$p$-completions.
\end{cor:hyperbolic}

However, the last statement of this Corollary could have been observed using Bousfield's result above.

These last Corollaries are interesting in light of M. Lackenby's result that

\begin{thm}\label{thm:lackenby1}~\cite[Theorem 1.12]{La3} Let $G$ be a finitely presented group. If the successive quotients of the derived $p$ series have $\textbf{linear growth rate}$ then $G$ is large, that is $G$ has a finite index subgroup that has a non-abelian free group as quotient.
\end{thm}

Here $\textbf{linear growth rate}$ means
$$
inf(dim_{\mathbb{Z}_p}(G^{(n)}/G^{(n+1)})/[G:G^{(n)}])>0.
$$

In 1975 William Dwyer extended Stallings' lower-central series theorem for $R=\mathbb{Z}$ by weakening the hypothesis on $H_2$ and indeed found precise conditions for when a \emph{specific} lower central series quotient was preserved ~\cite{Dw}. This was extended to the $R=\mathbb{Q}$ case by the authors ~\cite{CH2}. To motivate, for topologists, the philosophy of Dwyer's extension, consider that, for \emph{spaces}, the hypothesis that 
$$
H_2(\pi_1(Y);\mathbb{Z})\to H_2(\pi_1(X);\mathbb{Z})
$$
is surjective is equivalent to saying that $H_2(X;\mathbb{Z})$ is generated by classes from $H_2(Y;\mathbb{Z})$ together with \emph{spherical} classes, that is classes that can be represented by $2$-spheres. Roughly speaking, Dwyer showed that if 
``represented by $2$-spheres'' were replaced by ``represented by certain more general class of surfaces'', depending on a parameter $m$, then Stallings' Theorem still held for values of $n$ roughly up to the fixed value $m$. Refer to ~\cite{CH2} for a general discussion. For this purpose, for any group $B$, he defined an important subgroup 
$$
\Phi_m(B)\equiv \text{kernel}(H_2(B;\mathbb{Z})\to H_2(B/B_m;\mathbb{Z}))\subset H_2(B;\mathbb{Z}).
$$
Dwyer's ``filtration'' of $H_2(B)$ has equivalent more geometric formulations in terms of \emph{gropes} and played a crucial role in Freedman and Teichner's work on $4$-manifold topology that strengthened the results of Freedman-Quinn ~\cite{FT2}~\cite{FQ}~\cite{Kr1}~\cite{Kr2}~\cite{KT}. In actuality, Dwyer did not prove what we refer to as Dwyer's extensions of Stallings' theorems except in the case $R=\mathbb{Z}$. In particular, for the case $R=\mathbb{Z}_p$, he instead proved a version using group \textbf{cohomology} and defined a filtration of cohomology using Massey products. Therefore, to complete our analogy and to fill this historical gap, we include in Section~\ref{sec:Dwyerlower} a statement and proof of a direct (\textbf{homological}) extension of Stallings' theorem for the lower-central $p$-series, that follows quite easily from mimicking Dwyer's proof for the $R=\mathbb{Z}$ case.

We then define a new filtration $\Phi^{(m)}(B)$ of $H_2(B;\mathbb{Z}_p)$ that is appropriate to the derived $p$-series:
$$
\Phi^{(m)}(B)\equiv \text{image}(H_2(B^{(m)};\mathbb{Z}_p)\to H_2(B;\mathbb{Z}_p)).
$$
The more obvious analogue:
$$
\text{kernel}(H_2(B;\mathbb{Z}_p)\to H_2(B/B^{(m)};\mathbb{Z}_p))
$$
turns out to be incorrect. With this we can state our main theorem:

\newtheorem*{thm:dwyermain}{Theorem~\ref{thm:dwyermain}}
\begin{thm:dwyermain} Let $A$ be a finitely-generated group and $B$ be a finitely presented group. If $\phi:A\to B$ induces an isomorphism (respectively monomorphism) on $H_1(- ;\mathbb{Z}_p)$ and an
epimorphism $H_2(A;\mathbb{Z}_p)\to H_2(B;\mathbb{Z}_p)/\Phi^{(m)}(B)$, then for any $n\leq m+1$, $\phi$ induces an isomorphism (respectively monomorphism) $A/A^{(n)}\to B/B^{(n)}$. 
\end{thm:dwyermain}

The following corollary is a generalization of Stallings' ~\cite[Theorem 6.5]{St}. It should be compared to Dwyer's result ~\cite[Proposition 4.3]{Dw} whose hypothesis is in terms of Massey products and whose conclusion is in terms of the \textbf{restricted mod p lower central series}.

\newtheorem*{cor:independence}{Corollary~\ref{cor:independence}}
\begin{cor:independence} Let $B$ be a finitely presented group and $p$ a prime such that
$$
H_2(B^{(n-1)};\mathbb{Z})\to H_2(B;\mathbb{Z}_p)
$$
is surjective ($H_2(B;\mathbb{Z}_p)$ is generated by $B^{(n-1)}$-surfaces). Let $\{x_i\}$ be a finite set of elements of $B$ which is linearly independent in $H_1(B;\mathbb{Z}_p)$. Then the subgroup $A$ generated by $\{x_i\}$ has the free $p$-solvable group $F/F^{(n)}$ as quotient (where $F$ is free on $\{x_i\}$.
\end{cor:independence}
\begin{proof} [Proof of Corollary~\ref{cor:independence}] Consider the induced map $\phi:F\to B$ and observe that it satisfies the hypotheses of Theorem~\ref{thm:dwyermain}. Thus 
$$
F/F{(n)}\to B/B^{(n)}
$$
is injective and factors through $A/A^{(n)}$. Thus $F/F^{(n)}\cong A/A^{(n)}$.
\end{proof}

Also, in Section~\ref{sec:Topappl} we extend Stallings' application to link concordance to a more general equivalence relations on links.

\section{Basics}\label{sec:basics}

We collect a few elementary observations that we need in the next section. No originality is claimed.

The derived $p$-series subgroups are verbal hence functorial and in particular fully-invariant ~\cite[p.3]{St2}. It follows directly from the definition that $A^{\sss (n)}/A^{\sss (n+1)}$ is an abelian group every element of which has order dividing $p$. Next we note that the derived $p$-series has an equivalent formulation.

\begin{lem}\label{lem:homological1}
$A^{(n+1)}=$ ker$(A^{(n)}\overset{\pi\otimes 1}\longrightarrow (A^{(n)}/[A^{(n)},A^{(n)}])\otimes_\mathbb{Z}\mathbb{Z}_p)$.
\end{lem}
\begin{proof} Clearly
$$
A^{(n+1)}=[A^{\sss (n)},A^{\sss (n)}](A^{\sss (n)})^p\subset ker(\pi\otimes 1).
$$
On the other hand, if we tensor the exact sequence
$$
0\to \mathbb{Z}\overset{\cdot~p}\longrightarrow \mathbb{Z}\to \mathbb{Z}
_p\to 0
$$
with $C=A^{(n)}/[A^{\sss (n)},A^{\sss (n)}]$ we get the exact sequence
$$
C\overset{\cdot~p}\longrightarrow C\overset{id\otimes 1}\longrightarrow C\otimes \mathbb{Z}_p\to 0, 
$$
so the kernel of $id\otimes 1$ is $C^p$. Then, since $\pi\otimes 1$ factors as 
$$
A^{(n)}\overset{\pi}\to A^{(n)}/[A^{(n)},A^{(n)}]\overset{id\otimes 1}\longrightarrow A^{(n)}/[A^{(n)},A^{(n)}]\otimes_\mathbb{Z}\mathbb{Z}_p
$$
we see that the kernel of $\pi\otimes 1$is $\pi^{-1}(C^p)$, which is clearly $[A^{\sss (n)},A^{\sss (n)}](A^{\sss (n)})^p$. 
\end{proof}
\begin{lem}\label{lem:homological2} If $A$ is finitely generated, $A/A^{(n)}$ is a finite $p$-group each of whose elements has order dividing $p^{n}$.
\end{lem}
\begin{proof} Proceed by induction on $n$. If $n=0$, $A/A^{(0)}=1$ so the result holds. Suppose $A/A^{(n)}$ is a finite $p$-group each of whose elements has order dividing $p^{n}$. Then $A^{(n)}$ is a subgroup of finite index of the finitely generated group $A$, hence is itself finitely generated. Thus $A^{\sss (n)}/A^{\sss (n+1)}$ is a \emph{finitely generated} abelian group each of whose elements has order dividing $p$, and hence is finite. Then from the exact sequence below we see that $A/A^{(n+1)}$ is an extension of a finite $p$-group by a finite $p$-group, hence is itself a finite $p$-group and every element has order dividing $p^{n+1}$.

\begin{equation*}
\begin{CD}
1      @>>>    \f{\an}{\anp}  @>>>   \f A{\anp}  @>>>
\f A{\an}    @>>> 1\\
\end{CD}
\end{equation*}
\end{proof}

A crucial observation is that the derived series is related to homology as follows:

\begin{lem}\label{lem:homological3}
 $\f{\an}{\anp}\cong H_1(A;\mathbb{Z}_p[A/A^{(n)}])$ as right $\mathbb{Z}_p[A/A^{(n)}]$-modules.
\end{lem}

\begin{proof} This can be seen several ways. First consider the right-hand side of the equivalence. Note that the map $A\to A/A^{(n)}$ endows $\mathbb{Z}_p[A/A^{(n)}]$ with the structure of a $(\mathbb{Z}A,\mathbb{Z}_p[A/A^{(n)}])$-bimodule; similarly for $\mathbb{Z}[A/A^{(n)}]$. Thus
$$
H_1(A;\mathbb{Z}_p[A/A^{(n)}]) ~\text{and}~ H_1(A;\mathbb{Z}[A/A^{(n)}])
$$
are defined as right modules over $\mathbb{Z}_p[A/A^{(n)}]$ and $\mathbb{Z}[A/A^{(n)}]$ respectively. These modules have well-known interpretations as
$$
H_1(A^{(n)};\mathbb{Z}_p) ~\text{and}~ H_1(A^{(n)};\mathbb{Z})
$$
respectively. A topologist usually thinks of the group $A$ as being replaced by the Eilenberg-Maclane space $K(A,1)$. From this perspective the homology groups with twisted coefficients
$$
H_1(A;\mathbb{Z}_p[A/A^{(n)}]) ~\text{and}~ H_1(A;\mathbb{Z}[A/A^{(n)}])
$$
can be interpreted as merely the homology groups of the covering space, $K(A^{(n)},1)$, of $K(A,1)$ with fundamental group $A^{(n)}$, 
$$
H_1(A^{(n)};\mathbb{Z}_p) ~\text{and}~ H_1(A^{(n)};\mathbb{Z}).
$$
Here the group of covering translations is $A/A^{(n)}$, enabling us to view these as right modules over $\mathbb{Z}_p[A/A^{(n)}]$ or $\mathbb{Z}[A/A^{(n)}]$ respectively. Since the first homology is merely the abelianization, these homology modules have a purely group-theoretic description as 
$$
A^{(n)}/[A^{(n)},A^{(n)}]\otimes_\mathbb{Z}\mathbb{Z}_p ~\text{and}~ A^{(n)}/[A^{(n)},A^{(n)}]
$$
respectively. 
On the other hand, starting from the purely group-theoretic stand-point, $A^{(n)}/[A^{(n)},A^{(n)}]$ is an abelian group on which $A$ acts by conjugation ($x\to a^{-1}xa$) and $A^{(n)}$ acts trivially. Thus $A^{(n)}/[A^{(n)},A^{(n)}]$ has the structure of a  right $\mathbb{Z}[A/A^{(n)}]$-module. Similarly $A^{(n)}/[A^{(n)},A^{(n)}]\otimes_\mathbb{Z}\mathbb{Z}_p$ is a $\mathbb{Z}_p$-vector space on which $A$ acts by conjugation giving it the structure of a right $\mathbb{Z}_p[A/A^{(n)}]$-module. The point is that it is well known that these module structures are the same as those derived just above.

Using the algebraist's view of group homology, the fact that
$$
\an/[\an,\an]\cong H_1(A;\mathbb{Z}[A/\an])
$$
is a consequence of the definition $H_1(A;\mathbb{Z}[A/\an])\equiv Tor_1^{A}(\mathbb{Z}[A/\an],\mathbb{Z})$ and the easy observation that the latter is $Tor_1^{\an}(\mathbb{Z},\mathbb{Z})\cong \an/[\an,\an]$ \cite[Lemma 6.2]{HS}.

Now, since $\pi\otimes 1$ above is surjective with kernel $A^{(n+1)}$, it induces a group isomorphism
$$
A^{(n)}/A^{(n+1)}\overset{\pi\otimes 1}\longrightarrow A^{(n)}/[A^{(n)},A^{(n)}]\otimes_\mathbb{Z}\mathbb{Z}_p
$$
which endows $A^{(n)}/A^{(n+1)}$ with the structure of a right $\mathbb{Z}_p[A/A^{(n)}]$-module. These remarks complete the proof.
\end{proof}

\section{Dwyer's Theorem for the $p$-lower central series}\label{sec:Dwyerlower}

First, for historical completeness and to complete the analogy, we state and prove an improvement of Stallings' result for the \textbf{$p$-lower central series} that is suggested by Dwyer's work (although neither stated nor proved in this form by Dwyer). Dwyer proved a version for cohomology that employs a different $p$-series and uses a filtration defined in terms of Massey products~\cite[Theorem 3.1]{Dw}. It would be interesting to know how his result compares to that below. 

Dwyer's work suggests the following \emph{filtration} of $H_2(B;\mathbb{Z}_p)$. 

\begin{defn}\label{defn:Dwyerfiltration} Let $\Phi_{p,m}(B)$ denote the kernel of
$$
H_2(B;\mathbb{Z}_p)\to H_2(B/B_{p,m};\mathbb{Z}_p).
$$
Note that, if $n\leq m$ then $\Phi_{p,m}(B)\subset\Phi_{p,n}(B)$.
\end{defn}

\begin{thm}\label{thm:dwyerlower} If $\phi:A\to B$ induces an isomorphism on $H_1(- ;\mathbb{Z}_p)$ and an
epimorphism $H_2(A;\mathbb{Z}_p)\to H_2(B;\mathbb{Z}_p)/<\Phi_{p,m}(B)>$, then for any $n\leq m+1$,
$\phi$ induces an isomorphism $A/A_{p,n}\to B/B_{p,n}$. 
\end{thm}
\begin{proof}[Proof of Theorem~\ref{thm:dwyerlower}] We view $m$ as fixed and proceed by ``induction'' on
$n$. The case $n=2$ is clear since $A/A_{p,2}$ is merely
$H_1(A;\mathbb{Z}_p)$  and by hypothesis $\phi$ induces
an isomorphism on $H_1(-;\mathbb{Z}_p)$.
Now assume that the theorem holds for $n\leq m$, i.e. $\phi$ induces an isomorphism $A/\an \cong
B/\bn$. We will prove that it
holds for $n+1$. Since the $p$-lower central series is fully invariant, the diagram below exists and is commutative. 
\begin{equation*}
\begin{CD}
1      @>>>    \f{\aln}{\alnp}  @>>>   \f A{\alnp}  @>>>
\f A{\aln}    @>>> 1\\
&&     @VV \phi_n V        @VV\phi_{n+1}V       @VV\phi V\\
1      @>>>    \f{\bln}{\blnp}  @>>>   \f B{\blnp}  @>>>
\f B{\bln}    @>>> 1\\
\end{CD}
\end{equation*}
In light of the Five Lemma, it suffices to show that
$\phi_n$ is an isomorphism.

Now consider Stallings' exact sequence \cite[Theorem 2.1]{St}, where all homology groups have $\mathbb{Z}_p$-coefficients.
\begin{equation*}
\begin{CD}
 @>>>  H_2(A) @>\pi_A>> H_2(A/\aln)  @>\partial_A>>
\aln/\alnp   @>>> 0\\
&&    @VV \phi_* V        @VV(\phi_{n})_*V       @VV\phi_{n}V\\
 @>>>   H_2(B)  @>\pi_B>>    H_2(B/\bln)  @>\partial_B>>
\bln/\blnp    @>>> 0\\
\end{CD}
\end{equation*}
By our inductive hypothesis $(\phi_n)_*$ is an isomorphism. It follows immediately that $\phi_n$ is surjective. Suppose $\phi_n(a)=0$. Choose $\alpha$ such that $\partial_A(\alpha)=a$. Then $(\phi_{n})_*(\alpha)=\beta$ where $\partial_B(\beta)=0$ so $\beta= \pi_B(\gamma)$. Since, $n\leq m$, by our hypothesis $\gamma=\phi_*(\delta)+\epsilon$ for some $\delta\in H_2(A)$ and $\epsilon$ in the kernel of $\pi_B$. It follows that $\pi_A(\delta)=\alpha$. Thus $a=\partial_A(\pi_A(\delta))=0$, so $\phi_n$ is injective.
\end{proof}

\section{The Derived p-series Version of Stallings' and Dwyer's Theorems}\label{sec:Dwyer}

Now we move on to the derived $p$-series and our main results. First we need to define the correct analogue, for the derived $p$-series, of Dwyer's filtration of the second homology.

\begin{defn}\label{defn:Nsurfaces} Suppose $N$ is a normal subgroup of a group $B$. Let $\Phi^N(B)$ be the image of the inclusion-induced $H_2(N;\mathbb{Z}_p)\to H_2(B;\mathbb{Z}_p)$. Specifically if $N=B^{(m)}$ we abbreviate $\Phi^N(B)$ by $\Phi^{(m)}(B)$. Equivalently $\Phi^N(B)$ consists of those classes represented by what we call \textbf{$N$-surfaces of $B$}. An $N$-surface of $B$ is a continuous map $f:\Sigma\to K(B,1)$ of a compact oriented surface $\Sigma$ where $f_{|\partial \Sigma}$ factors through the standard map $z\to z^p$ on each circle (so $\Sigma$ represents a mod $p$ $2$-cycle) and such that $f_*(\pi_1(\Sigma))\subset N$.
\end{defn}

\begin{thm}\label{thm:dwyermain} Let $A$ be a finitely-generated group and $B$ be a finitely presented group. If $\phi:A\to B$ induces an isomorphism (respectively monomorphism) on $H_1(- ;\mathbb{Z}_p)$ and an
epimorphism $H_2(A;\mathbb{Z}_p)\to H_2(B;\mathbb{Z}_p)/<\Phi^{(m)}(B)>$, then for any $n\leq m+1$,
$\phi$ induces an isomorphism (respectively monomorphism) $A/A^{(n)}\to B/B^{(n)}$. 
\end{thm}

The following corollary is a precise analogue, for the derived $p$-series, of Stallings' theorem for the $p$-lower central series. It can also be proved directly from Stallings' theorem and in that proof it is seen that $B$ need only be finitely generated. A special case of this Corollary was shown in ~\cite[Theorem 1]{CO}.

\begin{cor}\label{cor:main} Let $A$ and $B$ be a finitely-generated groups. If $\phi:A\to B$ induces an isomorphism (respectively monomorphism) on $H_1(- ;\mathbb{Z}_p)$ and an
epimorphism on $H_2(- ;\mathbb{Z}_p)$, then for each finite $n$,
it induces an isomorphism (respectively monomorphism) $A/\an\to B/\bn$. $\phi$ also induces a monomorphism $A/A^{(\omega)} \subset B/B^{(\omega)}$. 
\end{cor}

\begin{proof}[Proof of Corollary \ref{cor:main}] In case $B$ is finitely presented, this is obviously follows from Theorem~\ref{thm:dwyermain}. This proof is independent of Stallings theorem. What follows is a more direct proof that relies on Stallings theorem.

\textbf{Alternative Proof of Corollary \ref{cor:main}:}~ We consider only the case that $\phi$ induces an isomorphism on $H_1(-;\mathbb{Z}_p)$. The monomorphism part follows by the same trick as used later in the proof of Theorem~\ref{thm:dwyermain}. By Lemma~\ref{lem:homological2} $A/\an$ and $B/\bn$ are finite $p$-groups. It is well known that any finite $p$-group is nilpotent. Moreover, there is a $k>0$ such that $(A/\an)_{p,k}=(B/\bn)_{p,k}=1$ ~\cite[Lemma 4.2]{St}. Thus, since the $p$-lower central series is fully invariant, the map $A\to A/\an$ factors through 
$A/A_{p,k}$ and so $A_{p,k}\subset \an$ and $B_{p,k}\subset \bn$. In examining the diagram below, Stallings' theorem implies that $\phi^k_*$ is an isomorphism.
\begin{equation*}
\begin{CD}
\f A{A_{p,k}}  @>{\phi^k_*}>>
\f B{B_{p,k}} \\
    @VV \pi_A V        @VV\pi_{B}V    \\
 \f A{\an}  @>{\phi}>>
\f B{\bn} \\
\end{CD}
\end{equation*}
Since the derived $p$-series is characteristic, $\phi$ induces an isomorphism
$$
\f {A/A_{p,k}}{(A/A_{p,k})^{(n)}}\overset{\phi}{\lra} \f {B/B_{p,k}}{(B/B_{p,k})^{(n)}}.
$$
It only remains to verify that
$$
\f A{\an}\cong \f {A/A_{p,k}}{(A/A_{p,k})^{(n)}}
$$
which follows from the lemma below.
\begin{lem}\label{lem:easygroup} If the kernel of the epimorphism $A\overset{f}{\lra}C$ is contained in $\an$ then $$
f^{(n)}:\f A{\an}\to \f C{C^{(n)}}
$$
is an isomorphism.
\end{lem}
\begin{proof}[Proof of Lemma~\ref{lem:easygroup}] Clearly $f^{(n)}$ is surjective. Suppose $f^{(n)}([a])=0$ so $f(a)\in C^{(n)}$. Since $f$ is surjective and the derived $p$-series is verbal, there is some $a'\in A^{(n)}$ such that $f(a')=f(a)$. Hence $a(a')^{-1}$ is in the kernel of $f$. Thus $a\in A{(n)}$.
\end{proof} 
\end{proof}

The proofs of the following corollaries are standard (compare ~\cite{St}).

\begin{cor}\label{cor:homologyequiv} Suppose $Y$ and $X$ are path-connected $CW$-complexes with $\pi_1(Y)$ and $\pi_1(X)$ finitely generated. If $f:Y\to X$ is a continuous map that induces
an isomorphism (respectively monomorphism) on $H_1(- ;\mathbb{Z}_p)$ and an
epimorphism on $H_2(- ;\mathbb{Z}_p)$, then for each finite $n$,
it induces an isomorphism (respectively monomorphism)
$$
f_*:\pi_1(Y)/\pi_1(Y)^{(n)}\to \pi_1(X)/\pi_1(X)^{(n)}.
$$
\end{cor}
\begin{proof}[Proof of Corollary~\ref{cor:homologyequiv}] Let $A=\pi_1(Y)$ and $B=\pi_1(X)$. Of course
$$
H_1(Y;\mathbb{Z}_p)\cong H_1(A;\mathbb{Z}_p).
$$
Since a $K(A,1)$ can be obtained from $Y$ by adjoining only cells of dimension $3$ and greater, the natural map
$$
H_2(Y;\mathbb{Z}_p)\to H_2(\pi_1(Y);\mathbb{Z}_p)
$$
is surjective. Thus the hypothesis that 
$$
f_*:H_2(Y;\mathbb{Z}_p)\to H_2(X;\mathbb{Z}_p)
$$
is surjective implies that
$$
f_*:H_2(A;\mathbb{Z}_p)\to H_2(B;\mathbb{Z}_p).
$$
Now apply Corollary~\ref{cor:main}.
\end{proof}

\begin{cor}\label{cor:homologycob} If $M$ and $N$ are compact $\mathbb{Z}_p$-homology cobordant manifolds then for each $n$
$$
\pi_1(M)/\pi_1(M)^{(n)}\cong \pi_1(N)/\pi_1(N)^{(n)}.
$$
\end{cor}
\begin{proof}[Proof of Corollary~\ref{cor:homologycob}] Let $C$ be the cobordism. Then the two inclusion maps are $\mathbb{Z}_p$ induce homology isomorphisms so Corollary~\ref{cor:homologyequiv} and be applied to each. Thus
$$
\pi_1(M)/\pi_1(M)^{(n)}\cong \pi_1(C)/\pi_1(C)^{(n)}\cong \pi_1(N)/\pi_1(N)^{(n)}.
$$
\end{proof}

\begin{cor}\label{cor:hyperbolic} For any closed orientable $3$-manifold $N$,
there exists a hyperbolic $3$-manifold $M$ such that for each $n$
$$
\pi_1(M)/\pi_1(M)^{(n)}\cong \pi_1(N)/\pi_1(N)^{(n)}
$$
and
$$
\pi_1(M)^{(n)}/\pi_1(M)^{(n+1)}\cong \pi_1(N)^{(n)}/\pi_1(N)^{(n+1)}.
$$
Hence the growth rate of the derived $p$-series of an arbitrary $3$-manifold is the same as that of a hyperbolic $3$-manifold with the same integral homology. Moreover $\pi_1(M)$ and $\pi_1(N)$ have isomorphic pro-$p$-completions.
\end{cor}

\begin{proof}[Proof of Corollary~\ref{cor:hyperbolic}] Most of this is is an immediate consequence of Corollary~\ref{cor:homologycob} and the fact, due to D. Ruberman that any closed orientable $3$-manifold is homology cobordant to a hyperbolic $3$-manifold ~\cite[Theorem 2.6]{Ru1}. However, if one wants the degree $1$ map, it is easier to apply the later result of Kawauchi : given $N$, there is a hyperbolic $M$ and a degree one map $f:M\to N$ that induces an isomorphism on all homology groups ~\cite{Ka1}\cite{Ka2} (rediscovered in ~\cite{BW}). Then simply apply Corollary~\ref{cor:homologyequiv}. The equivalence of the pro-$p$-completions follows immediately from Bousfield's result Corollary~\ref{cor:prop}.
\end{proof} 

The following is a generalization of Stallings' ~\cite[Theorem 6.5]{St}. It should be compared to Dwyer's result ~\cite[Proposition 4.3]{Dw} whose hypothesis is in terms of Massey products and whose conclusion is in terms of the \textbf{restricted mod p lower central series}.

\begin{cor}\label{cor:independence} Let $B$ be a finitely presented group and $p$ a prime such that
$$
H_2(B^{(n-1)};\mathbb{Z}_p)\to H_2(B;\mathbb{Z}_p)
$$
is surjective ($H_2(B;\mathbb{Z}_p)$ is generated by $B^{(n-1)}$-surfaces). Let $\{x_i\}$ be a finite set of elements of $B$ which is linearly independent in $H_1(B;\mathbb{Z}_p)$. Then the subgroup $A$ generated by $\{x_i\}$ has the free $p$-solvable group $F/F^{(n)}$ as quotient (where $F$ is free on $\{x_i\}$.
\end{cor}
\begin{proof} [Proof of Corollary~\ref{cor:independence}] Consider the induced map $\phi:F\to B$ and observe that it satisfies the hypotheses of Theorem~\ref{thm:dwyermain}. Thus 
$$
F/F{(n)}\to B/B^{(n)}
$$
is injective and factors through $A/A^{(n)}$. Thus $F/F^{(n)}\cong A/A^{(n)}$.
\end{proof}

\begin{proof}[Proof of Theorem \ref{thm:dwyermain}] We first consider the case that $\phi$ induces
an isomorphism on $H_1(-;\mathbb{Z}_p)$. We view $m$ as fixed and proceed by ``induction'' on
$n$. The case $n=1$ is clear since $A/A^{\sss (1)}$ is merely
$H_1(A;\mathbb{Z}_p)$ by Lemma~\ref{lem:homological3} and by hypothesis $\phi$ induces
an isomorphism on $H_1(-;\mathbb{Z}_p)$.
Now assume that the first claim holds for $n\leq m$, i.e. $\phi$ induces an isomorphism $A/\an \cong
B/\bn$. We will prove that it
holds for $n+1$.

Since the derived $p$-series is fully invariant, $\phi(A\np)\subset \bnp$. Hence the diagram below exists and is commutative. In light of the Five Lemma, it suffices to show that
$\phi$ induces an isomorphism $\an/\anp\to\bn/\bnp$.
\begin{equation*}
\begin{CD}
1      @>>>    \f{\an}{\anp}  @>>>   \f A{\anp}  @>>>
\f A{\an}    @>>> 1\\
&&     @VV \phi_n V        @VV\phi_{n+1}V       @VV\phi V\\
1      @>>>    \f{\bn}{\bnp}  @>>>   \f B{\bnp}  @>>>
\f B{\bn}    @>>> 1\\
\end{CD}
\end{equation*}
Now Lemma~\ref{lem:homological3} suggests that this problem can be translated into a homological one. In preparation, note that the inductive hypothesis is that $\phi$ induces an isomorphism $A/\an\to B/\bn$ and hence a ring isomorphism
$\phi:\mathbb{Z}_p[ A/\an]\to\mathbb{Z}_p[ B/\bn]$. Thus any left (right) $\mathbb{Z}_p[ B/\bn]$-module inherits a left (right) $\mathbb{Z}_p[ A/\an]$-module structure via $\phi$. The map $\phi$ also endows $\mathbb{Z}_p[ B/\bn]$ with the structure of a ($\mathbb{Z}_p[ A/\an]-\mathbb{Z}_p[ B/\bn]$)-bimodule with respect to which $\phi$ is a map of ($\mathbb{Z}_p[ A/\an]-\mathbb{Z}_p[ A/\an]$)-bimodules.

Now consider the following commutative diagram where the horizontal equivalences follow from Lemma~\ref{lem:homological3}.
\begin{equation*}
\begin{CD}
H_1(A;\mathbb{Z}_p[A/\an])  @>\cong>>  \f{\an}{\anp}\\ 
@VV \id\ox\phi V \\
H_1(A;\mathbb{Z}_p[B/\bn])\\
@VV \phi\ox\id V \\
H_1(B;\mathbb{Z}_p[B/\bn])  @>\cong>>  \f{\bn}{\bnp}\\ 
\\
\end{CD}
\end{equation*}
We claim that the composition
 $\phi_n=(\phi\ox\id)\circ(\id\ox\phi)$ is an isomorphism. Since $\phi:A/\an\to B/\bn$ is an isomorphism, $(\id\ox\phi)$ is clearly an isomorphism. Finally, we show that $\phi\ox\id$ is an isomorphism using Proposition~\ref{prop:2-connected} below (setting $\G = B/\bn, N=\bn$). Note that $B/\bn$ is a finite $p$-group by Lemma~\ref{lem:homological2}. This completes the proof of the ``isomorphism'' case of Theorem~\ref{thm:dwyermain}, modulo the proof of Proposition~\ref{prop:2-connected}.

\begin{prop}\label{prop:2-connected} Suppose $A$ is finitely-generated and $B$ finitely-related. Suppose $\phi:A\to B$ induces a monomorphism (respectively an isomorphism) on $H_1(-;\mathbb{Z}_p)$. Consider the coefficient system $\psi:B\to\G$ where $\G$ is a finite $p$-group. Suppose $N\subset\ker\psi$ and suppose that 
$$
H_2(A;\mathbb{Z}_p)\to H_2(B;\mathbb{Z}_p)/\Phi^N(B)
$$
is surjective. Equivalently, suppose that $H_2(B;\mathbb{Z}_p)$ is spanned by $\phi_*(H_2(A;\mathbb{Z}_p))$ together with a collection of $N$-surfaces. Then $\phi$ induces a monomorphism (respectively, an isomorphism)
$$
\phi_*: H_1(A;\mathbb{Z}_p\G)\lra H_1(B;\mathbb{Z}_p\G).
$$
\end{prop}

Before, proving Proposition~\ref{prop:2-connected}, we finish the proof of Theorem~\ref{thm:dwyermain}. Suppose we assume only that $\phi$ induces a \emph{monomorphism} on $H_1(- ;\mathbb{Z}_p)$. We are grateful to Kent Orr for suggesting the idea for this argument. Since $H_1(B;\mathbb{Z}_p)$ is a $\mathbb{Z}_p$-vector space, it decomposes as image($\phi$)$\oplus C$ for some vector space $C$ with basis $\{c_i|i\in \mathcal{C}\}$. Let $F$ be the free group on $\mathcal{C}$ and $\tilde A=A*F$. We can extend $\phi$ to a map $\psi:\tilde A\to B$, by setting $\psi(x_i)=c_i$. Observe that $\psi$ induces an isomorphism on $H_1(- ;\mathbb{Z}_p)$. Moreover there is an obvious retraction $r:\tilde{A}\to A$, under which $H_2(\tilde{A};\mathbb{Z}_p)$ maps onto $H_2(A;\mathbb{Z}_p)$, so we recover our hypothesis on $H_2$ as well. Thus, by the first part of the Theorem, for any finite $n$,
$\psi$ induces an isomorphism $\tilde{A}/\tilde{\an}\to B/\bn$. Since the derived p-series is fully invariant, there is a retraction $r:\tilde{A}/\tilde{\an}\to A/\an$ showing that $A/\an\to \tilde{A}/\tilde{\an}$ is injective. Thus $\phi$ induces a monomorphism $A/A^{(n)} \to
B/B^{(n)}$. This completes the proof of Theorem~\ref{thm:dwyermain}, modulo the proof of Proposition~\ref{prop:2-connected}.

\begin{proof}[Proof of Proposition~\ref{prop:2-connected}] We first remark that since $\phi_*$ is certainly a map of $\mathbb{Z_p}\G$-modules, it suffices to show that $\phi_*$ induces a bijection of $\mathbb{Z}_p-\textbf{vector spaces}$. Then, since any finite $p$-group is poly-(cyclic of p-power order), it would suffice to prove the theorem for $\G=\mathbb{Z}_{p^r}$. The general proof is not much more difficult.

We offer a ``topological'' proof. We may consider that $A$, $B$ are Eilenberg-Maclane spaces with a finite number of 1-cells for $A$ and 2-cells for $B$ and that $\phi$ is cellular. By replacing the chosen $K(B,1)$ with the mapping cylinder of $\phi$, we may assume that $K(A,1)$ is a subcomplex of $K(B,1)$ and that the latter has a finite $2$-skeleton. Thus with such a cell structure the ordinary cellular chain complex $\ov\SC_*=\ov C_*(B,A;\mathbb{Z}_p)$ and the twisted relative cellular chain complex $\SC_*\equiv\SC_*(B,A;\mathbb{Z}_p\G)$ which will be finitely generated in dimension 2. If one thinks of $\psi$ as inducing a principal $\G$-bundle $B_{\G}$ over $B$ and $\psi\circ\phi$ as inducing one, $A_{\G}$, over $A$, then $\SC_*$ is merely the relative cellular chain complex for $(B_{\G},A_{\G})$ with $\mathbb{Z}_p\G$ coefficients. Let $\ov\SC_*=\SC_*\otimes_{\mathbb{Z}_p\G}\mathbb{Z}_p$ with $\pi_{\#}:\SC_*\to\ov\SC_*$ and note that $\ov\SC_*$ may be identified with the ordinary cellular chain complex of $(B,A)$ with $\mathbb{Z}_p$-coefficients.

Let $\{\ov\Sigma_s\mid s\in S\}$ denote the collection of $N$-surfaces in $B$ and let $\{\bar x_s\}$ denote the 2-cycles in $\ov C_2(B,A)$ represented by $\{\ov\Sigma_s\}$. Examining the exact sequence below 
$$
H_2(A;\mathbb{Z}_p)\overset{\phi_*}{\lra} H_2(B;\mathbb{Z}_p)\lra H_2(B,A;\mathbb{Z}_p) \overset{\p_*}{\lra} H_1(A;\mathbb{Z}_p)\overset{\phi_*}{\lra} H_1(B;\mathbb{Z}_p)
$$
where by hypothesis, $\phi_*$ is injective on $H_1(-;\mathbb{Z}_p)$, one sees that our hypotheses are engineered precisely so that $H_2(\ov\SC_*)$ is spanned by $\{[\bar x_s]\}$. Hence
$$
\rank_{\mathbb{Z}_p}(H_2(\ov\SC_*))/\<\{[\bar x_s]\}\>) = 0.
$$

Now since $N\sbq\ker\psi$, the $\ov\Sigma_s$ are also $(\text{ker}\psi)$-surfaces and so lift to $B_{\G}$. Choose lifts $\{\Sigma_s\mid s\in S\}$. That is, by definition, the $\{\ov\Sigma_s\}$  can be lifted to $H_2(N;\mathbb{Z}_p)$ and hence can be lifted to represent classes $\{[\Sigma_s]\}$ in $H_2(B_{\G};\mathbb{Z}_p)\cong H_2(B;\mathbb{Z}_p\G)$, and hence represent classes in $H_2(B,A;\mathbb{Z}_p\G)$. Let $\{x_s\mid s\in S\}$ denote the 2-cycles of $C_2(B,A;\mathbb{Z}_p\G)$ represented by $\{\Sigma_s\}$. Note that $\bar x_s=\pi_{\#}(x_s)$. Now consider the exact sequence
$$
H_2(B;\mathbb{Z}_p\G)\overset{\pi_*}{\lra} H_2(B,A;\mathbb{Z}_p\G)\overset{\p_*}{\lra} H_1(A;\mathbb{Z}_p\G)
\overset{\phi_*}{\lra} H_1(B;\mathbb{Z}_p\G).
$$
Our goal, namely proving that $\phi_*$ is injective, is now seen to be equivalent to showing $\pi_*$ is surjective. Let $\<\{[x_s]\}\>$ denote the $\mathbb{Z}_p\G$-submodule generated by $\{[x_s]\}$. Since $\{[x_s]\}\subset\image\pi_*$, it suffices to show that
$$
\rank_{\mathbb{Z}_p}(H_2(\SC_*))/\<\{[x_s]\}\>) = 0.
$$
Therefore we are reduced to proving a homological statement:
$$
\rank_{\mathbb{Z}_p}(H_2(\ov\SC_*))/\<\{[\bar x_s]\}\>) = 0 \Rightarrow \rank_{\mathbb{Z}_p}(H_2(\SC_*))/\<\{[x_s]\}\>) = 0.
$$
We can simplify this even more as follows. We can define a projective chain complex $\SD_*=\{D_q,d_q\}$ such that $H_2(\SD_*)\cong H_2(\SC_*)/\<[x_s]\mid s\in S\>$. Set $D_{3}=(\oplus_{s\in S}R\G)\bigoplus C_{3}$ and otherwise set $D_q=C_q$. Let $d_{3}:D_{3}\to D_2$ be defined by $d_{3}(e_s,y)=x_s+\p_{3}(y)$ where $\{e_s\}$ is a basis of $(R\G)^s$ and $y\in C_{3}$, and $d_{4}:D_{4}\to D_{3}$  by $d_{4}(z)=(0,\p_{4}(z))$ for $z\in D_{4}=C_{4}$. Then the $2$-cycles of $\SD_*$ are the same as those of $\SC_*$ while the group of $2$-boundaries is larger (includes $\{x_s\}$). Hence $H_2(\SD_*)\cong H_2(\SC_*)/\<[x_s]\mid s\in S\>$ as claimed. Note that $D_2=C_2$ is finitely generated. Similarly we can define a chain complex $\ov\SD_*$ which agrees with $\ov\SC_*$ except in dimension $3$ where $\ov D_{3}=(\oplus_{s\in S}R)\bigoplus\ov\SC_*$ with $\bar d_{3}:\ov D_{3}\to\ov D_2$ given by $d_{3}(\bar e_s,\bar y)=\bar x_s+\ov\p_{3}(\bar y)$ for $\{\bar e_s\}$ a basis of $R^s$ and $\bar y\in\ov C_{3}$. Then just as above:
$$
H_2(\ov\SD_*)\cong H_2(\ov\SC_*)/\<[x_s] ~| s\in S\>.
$$
With this translation we are reduced to proving:
$$
\rank_{\mathbb{Z}_p}(H_2(\ov\SD_*))) = 0 \Rightarrow \rank_{\mathbb{Z}_p}(H_2(\SD_*))) = 0.
$$
Moreover the chain map $\pi:\SC_*\to\ov\SC_*$ extends to $\tl\pi:\SD_*\to\ov\SD_*$ by setting $\tl\pi(e_s,y)=(\bar e_s,\pi(y))$, and one sees that $\ov\SD_*=\SD_*\otimes_{\mathbb{Z}_p\G}\mathbb{Z}_p$.
and then apply Corollary~\ref{cor:rank} to this chain complex. 

The desired result now follows immediately from Corollary~\ref{cor:rank} (proof postponed).

\begin{cor}\label{cor:rank}Suppose $\G$ is a finite $p$-group and $\SD_*$ is a projective  right $\mathbb{Z}_p\G$ chain complex with $D_q$ finitely generated. Then
$$
\rank_{\mathbb{Z}_p}H_q(\SD_*)\le |\G|\rank_{\mathbb{Z}_p} H_q(\SD_*\otimes_{\mathbb{Z}_p\G}\mathbb{Z}_p).
$$
\end{cor}

This concludes the proof that $\phi_*$ is \emph{injective} with $\mathbb{Z}_p\G$ coefficients, modulo the proof of Corollary~\ref{cor:rank}.

Before proving Corollary~\ref{cor:rank}, consider the case that, in addition, $\phi_*$ is an \emph{isomorphism} on $H_1(-;\mathbb{Z}_p)$. Since $B$ is finitely generated and we can assume that $K(A,1)$ has only a finite number of $0$ cells, after the mapping cylinder construction, it follows that $\SC_*(B,A;\mathbb{Z}_p\G)$ is finitely generated in dimension $1$. Again, by Corollary~\ref{cor:rank}, we see that
$$
\rank_{\mathbb{Z}_p}H_1(\SC_*)\leq |\G|\rank_{\mathbb{Z}_p}H_1(\ov \SC_*)
$$
and thus
$$
\rank_{\mathbb{Z}_p}H_1(B,A;\mathbb{Z}_p)\leq |\G|\rank_{\mathbb{Z}_p}H_1(B,A;\mathbb{Z}_p)=0.
$$

It follows that $\phi$ induces an epimorphism and hence an isomorphism
$$
\phi_*: H_1(A;\mathbb{Z}_p\G)\lra H_1(B;\mathbb{Z}_p\G).
$$

This completes the proof of Proposition~\ref{prop:2-connected} modulo the proof of Corollary~\ref{cor:rank}. In order to prove Corollary~\ref{cor:rank} we will need the following result of R. Strebel and a modest corollary.

\begin{prop}\label{injective} {(R. Strebel~\cite[Lemma 1.10]{Str})} Suppose $p$ is a prime integer and $\G$ is a
finite $p$-group. Any map between projective right
$\mathbb{Z}_p\G$-modules whose image under the functor
$-\otimes_{\mathbb{Z}_p\G}\mathbb{Z}_p$ is injective, is itself injective.
\end{prop}

Remarks: Strebel established this property for a much larger class of groups which Howie and Schneebli showed was precisely the class of all locally p-indicable groups ~\cite{HoS}.

\begin{lem}\label{lem:strebel} Suppose $\tl f:M\to N$ is a homomorphism
between projective $\mathbb{Z}_p\G$-modules with $\Gamma$  a finite p-group and let $f=\tl
f\ox\id$ be the induced homomorphism of $\mathbb{Z}_p$-vector spaces
$M\ox_{\mathbb{Z}_p\G}\mathbb{Z}_p\to N\ox_{\mathbb{Z}_p\G}\mathbb{Z}_p$. Then
$\rank_{\mathbb{Z}_p}(\image\tl f)\ge$ $|\G|\rank_{\mathbb{Z}_p}(\image f)$ where $|\G|$ is the order of $\G$.
\end{lem}

\begin{proof}[Proof of Lemma~\ref{lem:strebel}] By the rank of a homomorphism
we shall mean the rank of its image. Suppose that
$\rank_{\mathbb{Z}_p} f\geq r$. Then there is a {\it
monomorphism} $g:\bz_p^r\to N\ox_{\bz_p\G}\bz_p$ whose image is a
subgroup of $\image f$. If $e_i$, $1\le i\le r$ is a basis of
$\mathbb{Z}_p^r$, choose $M_i\in M\ox_{\mathbb{Z}_p\G}\mathbb{Z}_p$ such
that $f(M_i)=g(e_i)$. Since the ``augmentation'' $\e_M:M\to
M\ox_{\mathbb{Z}_p\G}\mathbb{Z}_p$ is surjective there exist elements
$m_i\in M$ such that $\e_M(m_i)=M_i$. Consider the map $\tl
g:(\mathbb{Z}_p\G)^r\to N$ defined by sending the $i^{\supth}$ basis
element to $\tl f(m_i)$. The augmentation of $\tl g$, $\tl
g\ox\id$, is the map
$(\mathbb{Z}_p\G)^r\ox_{\mathbb{Z}_p\G}\mathbb{Z}_p\to
N\ox_{\mathbb{Z}_p\G}\mathbb{Z}_p$ that sends $e_i$ to $\e_N(\tl f(m_i))=
f(\e_M(m_i))=g(e_i)$ and thus is seen to be
identifiable with $g$. In particular $\tl g\ox\id$ is a
monomorphism, and thus by Proposition~\ref{injective} ,
$\tl g$ is a monomorphism. Since the image of $\tl g$ lies in the
image of $\tl f$, $\tl g$ yields a monomorphism from
$(\mathbb{Z}_p\G)^r$ into the image of $\tl f$, showing that the
$\mathbb{Z}_p$-rank of $\image\tl f$ is at least $r$ times the order of $\G$.
\end{proof}

\begin{proof}[Proof of Corollary~\ref{cor:rank}] Let $\{\ov\SD_*\}=\{\ov D_q,\ov\p_q\}=\{D_q\otimes_{\mathbb{Z}_p\G}\mathbb{Z}_p,\p_q\ox\id\}$ , $\ov r_q= rank_{\mathbb{Z}_p}\ov D_q$ and let $r_q= rank_{\mathbb{Z}_p}D_q$. Since $D_q$ is finitely generated and projective and $\G$ is finite, we claim that $r_q=|\G|\ov r_q$. This is obvious if $D_q$ is free and requires a small argument if not (same as that in the proof of ~\cite[Corollary 2.8]{CH2}).

Now observe
\begin{align*}
\rank_{\mathbb{Z}_p}H_q(\SD_*)  &= \rank_{\mathbb{Z}_p}(\ker\p_q) - \rank_{\mathbb{Z}_p}(\image\p_{q+1})\\
&= r_q - \rank_{\mathbb{Z}_p}(\image\p_q) - \rank_{\mathbb{Z}_p}(\image\p_{q+1})\\
&\le |\G|\ov r_q - |\G|\rank_{\mathbb{Z}_p}(\image\ov\p_q) - |\G|\rank_{\mathbb{Z}_p}(\image\ov\p_{q+1})\\
&= |\G|(\rank_{\mathbb{Z}_p}(\ker\ov\p_q) - \rank_{\mathbb{Z}_p}(\image\ov\p_{q+1}))\\
&= |\G|\rank_{\mathbb{Z}_p} H_q(\ov\SD_*),
\end{align*}
where the inequality follows from two applications of Lemma~\ref{lem:strebel} and our observation that $r_q=|\G|\ov r_q$. This completes the proof of Corollary~\ref{cor:rank}.
\end{proof}

This finishes the proofs of Proposition~\ref{prop:2-connected} and hence that of Theorem~\ref{thm:dwyermain}.
\end{proof}
\end{proof}

The following is a useful consequence.

\begin{cor}\label{cor:n-conn} Suppose $\Gamma$ is a finite $p$-group. \ 
\begin{itemize}  
\item[a)] If $C_{\ast}$ is a non-negative right $\mathbb{Z}_p\Gamma$ chain complex that is
finitely-generated and projective in dimensions $0\le i\le m$ such that
$H_i(C_{\ast}\otimes_{\mathbb{Z}_p\G}\mathbb{Z}_p)=0$ for $0\le i\le m$, then
$H_i(C_{\ast})=0$ for $0\le i\le m$. 
\item[b)] If $f:Y\to X$ is a continuous map between connected CW complexes, X having a finite
$m$-skeleton and $Y$ a finite $(m-1)$-skeleton, which is $m$-connected on $\mathbb{Z}_p$ homology, and
$\phi:\pi_1(X)\to\G$ is a coefficient system, then $f$ is
$m$-connected on homology with $\mathbb{Z}_p\Gamma$-coefficients. In particular, if $Y_{\G}$ and $X_\G$ denote the induced $\G$-covering spaces $Y$ and $X$ then any lift $\tilde f: Y_{\G}\to X_\G$ of $f$ is $m$-connected on $\mathbb{Z}_p$ homology.
\end{itemize}
\end{cor}

\begin{proof}   Let $\e:\mathbb{Z}_p\G\to\mathbb{Z}_p$ be the augmentation and $\e(C_{\ast})$ denote
$C_{\ast}\otimes_{\mathbb{Z}_p\G}\mathbb{Z}_p$. Since $\e(C_{\ast})$ is acyclic up to dimension~$m$,
there is a $``$partial'' chain homotopy
$$
\{h_i:\e(C_{\ast})_i\to\e(C_{\ast})_{i+1}\mid0\le i\le m\}
$$ between the identity and the zero chain homomorphisms. By this we mean that
$\partial h_i+h_{i-1}\partial =\id$ for $0\le i\le m$.

Since $C_i \xrightarrow{\e} \e(C_i)$ is surjective and $C_i$ is projective, $h_i\circ \e:C_i\to \e(C_i)$ can be lifted to $\widetilde h_i:C_i\to C_i$, so that
$\e\circ\widetilde h_i=h_i\circ \e$. In this manner $h$ can be lifted to a partial chain homotopy
$\{\widetilde h_i\mid0\le i\le m\}$ on
$C_{\ast}$ between some partial chain map $\{f_i\mid0\le i\le m\}$ and the
zero map. Moreover $\e(f_i)$ is the identity map on $\e(C_{\ast})_i$, and in
particular, is injective. Thus, by Proposition~\ref{injective}, $f_i:C_i\to C_i$ is an
injective map of $\mathbb{Z}_p\Gamma$-modules for each $i$. Note $C_i$ is also a $\mathbb{Z}_p$-vector space. Since $C_i$ is a finitely generated projective $\mathbb{Z}_p\Gamma$-module, it is a direct summand of a finitely generated free $\mathbb{Z}_p\Gamma$-module. Since $\Gamma$ is a finite group, $\mathbb{Z}_p\Gamma$ is a finite dimensional $\mathbb{Z}_p$-vector space. It follows that $C_i$ is a \emph{finite dimensional} $\mathbb{Z}_p$-vector space. Any injective map between finite dimensional vector spaces of the same rank is necessarily bijective. Thus $f_i$ is bijective. Thus $f_i$ is a bijective morphism of $\mathbb{Z}_p\Gamma$-modules, hence an isomorphism of $\mathbb{Z}_p\Gamma$-modules for each $i$. Consequently, $\widetilde h_i$ is a partial
chain homotopy on
$C_{\ast}$ between the zero map and the partial chain map
$f_i$. To see that $C_{\ast}$ is acyclic, suppose $z\in C_i$ is a cycle, $0\leq i \leq m$. Then $f_i(w)=z$ for some $w\in C_i$. Thus
$$
\partial h_i(w)+ h_{i-1}\partial(w)=f_i(w)=z.
$$
If $w$ is a cycle then we see immediately that $z$ is a boundary. To see that $w$ is a cycle, apply boundary to both sides of the above equation to yield $\partial h_{i-1}\partial(w)=0$. Hence
$$
0=(\partial h_{i-1})\partial(w)=(f_{i-1}-h_{i-2}\partial)\partial(w)=f_{i-1}(\partial w)
$$
which implies that $\partial(w)=0$ since $f_{i-1}$ is an isomorphism.
Hence we have shown that
$C_{\ast}$ is acyclic up to and including dimension~$m$.

The second statement follows from applying this to the relative cellular chain
complex associated to the mapping cylinder of $f$.
\end{proof}
\section{Topological Applications}\label{sec:Topappl}

Stallings' and Dwyer's theorems have been instrumental in the study of homology cobordism of $3$-manifolds and in particular in the study of link concordance. Our current results may be similarly applied. By a \textbf{link} $L$ of $m$ components we mean an oriented, ordered collection of $m$ of circles disjointly and smoothly embedded in $S^3$. Two links $L_0$ and $L_1$ are \textbf{concordant} if there exist $m$ annuli disjointly embedded in $S^3\times [0,1]$, restricting to yield $L_j$ on $S^3\times \{j\}$, $j=0,1$. The complement of the union of annuli  is easily seen, by Alexander Duality, to be a product on integral homology, and thus the fundamental groups of the link complements are related by Stallings' theorems. Recently, several other weaker equivalence relations on knots and links have been considered and found to be useful in understanding knot and link concordance ~\cite{COT}~\cite{COT2}~\cite{CT}~\cite{ConT1}~\cite{ConT2}~\cite{CT}~\cite{KT}~\cite{Kr1}~\cite{Ha2}~\cite{T}. These equivalence relations involved replacing the annuli in the definition of concordance by surfaces equipped with some extra structure depending on a parameter. As the parameter increases, these surfaces are to be viewed as ``approximating'' annuli, hence yielding ``filtrations'' of link concordance. Below we show that our results generalize Stallings' results on link concordance to some of these more general equivalence relations. 

Recall that Stallings showed that concordant links have exteriors whose fundamental groups are isomorphic modulo any term of the lower central $p$-series. We can prove an analogue for the derived $p$-series and moreover use our Dwyer-type theorem to generalize his result to the following equivalence relation that is weaker than concordance.

\begin{defn}\label{ncobordant}The $m$-component links $L_0$ are \emph{$(n)$-$\mathbb{Z}_p$-cobordant} if there exist compact oriented surfaces $\Sigma_i, 1\leq i\leq m$, properly and disjointly embedded in $S^3\times [0,1]$, restricting to yield $L_j$ on $S^3\times \{j\}$, $j=0,1$, such that, for each $i$, for some set of circles $\{a_j,b_j\}$ representing a symplectic basis of curves for $\Sigma_i$, the image of each of the loops $\{a_j,b_j\}$ in $\pi_1((S^3\times [0,1])- \coprod \Sigma_i)\equiv \pi_1(E)$ is contained in $\pi_1(E)^{(n)}$ (use the unique `unlinked' normal vector field on $\Sigma_i$ to push off). 
\end{defn}

\begin{prop}\label{injectivity}If $L_0$ and $L_1$ are $(n)$-$\mathbb{Z}_p$-cobordant and $A=\pi_1(S^3-L_0)$, $\overline{A}=\pi_1(S^3-L_1)$, then
$$
A/A^{(n)}\cong \overline{A}/\overline{A}^{(n)}.
$$ 
\end{prop}

\begin{proof} We use the notation of Definition~\ref{ncobordant} and Proposition~\ref{injectivity}. Let $B=\pi_1(E)$. By hypothesis, for each $1\leq i\leq m$ there exist symplectic bases of circles $\{a_{ij},b_{ij}\}$ for $\Sigma_i$ whose push-offs $\{a^+_{ij},b^+_{ij}\}$ into $E$ lie in $\pi_1(E)^{(n)}= B^{(n)}$. The key observation is that $H_2(E;\mathbb{Z}_p)$ is generated by the tori $\{a^+_{ij}\times S^1_i, b^+_{ij}\times S^1_i\}$ where $S^1_i$ is a fiber of the normal circle bundle to $\Sigma_i$, together with the $m$ tori $L_0\times S^1_i$ that live in $S^3-L_0$. Thus the cokernel of the map $H_2(A;\mathbb{Z}_p)\to H_2(B;\mathbb{Z}_p)$ is generated by the former collections. Since $[a^+_{ij}]\in B^{(n)}$, $a^+_{ij}$ bounds a $B^{(n-1)}$-surface $S_{ij}$ mapped into $E$ (that is $\pi_1(S_{ij})\subset B^{(n-1)}$). If we cut open the torus $a^+_{ij}\times S^1_i$ along $a^+_{ij}$ and adjoin two oppositely oriented copies of $S_{ij}$, we obtain a (mapped in) surface that is homologous to $a^+_{ij}\times S^1_i$ and is also a $B^{(n-1)}$-surface. Similarly for the tori $b^+_{ij}\times S^1_i$. Therefore $A\to B$ satisfies the hypotheses of Theorem~\ref{thm:dwyermain} for $n-1$. Symmetrically, the same is true for $\overline{A}\to B$. The theorem follows immediately.
\end{proof}

\bibliographystyle{plain}
\bibliography{mybib3}
\end{document}

%% file: amsmath-defs.tex
\def\e{\epsilon}

\def\f{\frac}

\def\G{\Gamma}

\def\lb\{{\left\{}
\def\la{\lambda}
\def\La{\Lambda}
\def\lla{\longleftarrow}
\def\lm{\limits}
\def\lra{\longrightarrow}
\def\dllra{\Longleftrightarrow}
\def\llra{\longleftrightarrow}
\def\n{\nabla}
\def\ngth{\negthickspace}
\def\ola{\overleftarrow}
\def\Om{\Omega}
\def\om{\omega}
\def\op{\oplus}
\def\oper{\operatorname}
\def\oplm{\operatornamewithlimits}
\def\ora{\overrightarrow}
\def\ov{\overline}
\def\ova{\overarrow}
\def\ox{\otimes}
\def\p{\partial}
\def\rb\}{\right\}}
\def\s{\sigma}
\def\sbq{\subseteq}
\def\spq{\supseteq}
\def\sqp{\sqsupset}
\def\supth{{\text{th}}}
\def\T{\Theta}
\def\th{\theta}
\def\tl{\tilde}
\def\thra{\twoheadrightarrow}
\def\un{\underline}
\def\ups{\upsilon}
\def\vp{\varphi}
\def\wh{\widehat}
\def\wt{\widetilde}
\def\x{\times}
\def\z{\zeta}
\def\({\left(}
\def\){\right)}
\def\[{\left[}
\def\]{\right]}
\def\<{\left<}
\def\>{\right>}

\def\tec{Teichm\"uller\ }
\def\sconr{\hbox{\medspace\vrule width 0.4pt height 4.7pt depth
0.4pt \vrule width 5pt height 0pt depth 0.4pt\medspace}}

\def\SA{\mathcal A}
\def\SB{\mathcal B}
\def\SC{\mathcal C}
\def\SD{\mathcal D}
\def\SE{\mathcal E}
\def\SF{\mathcal F}
\def\SG{\mathcal G}
\def\SH{\mathcal H}
\def\SI{\mathcal I}
\def\SJ{\mathcal J}
\def\SK{\mathcal K}
\def\SL{\mathcal L}
\def\SM{\mathcal M}
\def\SN{\mathcal N}
\def\SO{\mathcal O}
\def\SP{\mathcal P}
\def\SQ{\mathcal Q}
\def\SR{\mathcal R}
\def\SS{\mathcal S}
\def\ST{\mathcal T}
\def\SU{\mathcal U}
\def\SV{\mathcal V}
\def\SW{\mathcal W}
\def\SX{\mathcal X}
\def\SY{\mathcal Y}
\def\SZ{\mathcal Z}

\newcommand{\BA}{\ensuremath{\mathbf A}}
\newcommand{\BB}{\ensuremath{\mathbf B}}
\newcommand{\BC}{\ensuremath{\mathbf C}}
\newcommand{\BD}{\ensuremath{\mathbf D}}
\newcommand{\BE}{\ensuremath{\mathbf E}}
\newcommand{\BF}{\ensuremath{\mathbf F}}
\newcommand{\BG}{\ensuremath{\mathbf G}}
\newcommand{\BH}{\ensuremath{\mathbf H}}
\newcommand{\BI}{\ensuremath{\mathbf I}}
\newcommand{\BJ}{\ensuremath{\mathbf J}}
\newcommand{\BK}{\ensuremath{\mathbf K}}
\newcommand{\BL}{\ensuremath{\mathbf L}}
\newcommand{\BM}{\ensuremath{\mathbf M}}
\newcommand{\BN}{\ensuremath{\mathbf N}}
\newcommand{\BO}{\ensuremath{\mathbf O}}
\newcommand{\BP}{\ensuremath{\mathbf P}}
\newcommand{\BQ}{\ensuremath{\mathbf Q}}
\newcommand{\BR}{\ensuremath{\mathbf R}}
\newcommand{\BS}{\ensuremath{\mathbf S}}
\newcommand{\BT}{\ensuremath{\mathbf T}}
\newcommand{\BU}{\ensuremath{\mathbf U}}
\newcommand{\BV}{\ensuremath{\mathbf V}}
\newcommand{\BW}{\ensuremath{\mathbf W}}
\newcommand{\BX}{\ensuremath{\mathbf X}}
\newcommand{\BY}{\ensuremath{\mathbf Y}}
\newcommand{\BZ}{\ensuremath{\mathbf Z}}

\def\bba{{\mathbb A}}
\def\bbb{{\mathbb B}}
\def\bbc{{\mathbb C}}
\def\bbd{{\mathbb D}}
\def\bbe{{\mathbb E}}
\def\bbf{{\mathbb F}}
\def\bbg{{\mathbb G}}
\def\bbh{{\mathbb H}}
\def\bbi{{\mathbb I}}
\def\bbj{{\mathbb J}}
\def\bbk{{\mathbb K}}
\def\bbl{{\mathbb L}}
\def\bbm{{\mathbb M}}
\def\bbn{{\mathbb N}}
\def\bbo{{\mathbb O}}
\def\bbp{{\mathbb P}}
\def\bbq{{\mathbb Q}}
\def\bbr{{\mathbb R}}
\def\bbs{{\mathbb S}}
\def\bbt{{\mathbb T}}
\def\bbu{{\mathbb U}}
\def\bbv{{\mathbb V}}
\def\bbw{{\mathbb W}}
\def\bbx{{\mathbb X}}
\def\bby{{\mathbb Y}}
\def\bbz{{\mathbb Z}}

%% file: pderivedseries4_1_07.bbl
\begin{thebibliography}{10}

\bibitem{BW}
Michel Boileau and Shicheng Wang.
\newblock Non-zero degree maps and surface bundles over {$S\sp 1$}.
\newblock {\em J. Differential Geom.}, 43(4):789--806, 1996.

\bibitem{B}
A.~Bousfield.
\newblock {\em Homological localization towers for groups and $\pi$-modules}.
\newblock Number 186 in Memoirs of the American Math. Soc. American Math. Soc.,
  Providence, Rhode Island, 1977.

\bibitem{CH2}
T.~Cochran and S.~Harvey.
\newblock Homology and derived series of groups ii: Dwyer's theorem.
\newblock preprint' http://xxx.lanl.gov/abs/math.GT/0609484.

\bibitem{COT}
T.~Cochran, K.~Orr, and P.~Teichner.
\newblock Knot concordance, whitney towers and $l^2$-signatures.
\newblock {\em Annals of Math.}, 157:433--519, 2003.

\bibitem{COT2}
T.~Cochran, K.~Orr, and P.~Teichner.
\newblock Structure in the classical knot concordance group.
\newblock {\em Comment. Math. Helv.}, pages 105--123, 2004.

\bibitem{CT}
T.~Cochran and P.~Teichner.
\newblock Knot concordance and von neumann $\rho$-invariants.
\newblock {\em Duke Math. Journal}.
\newblock to appear, mathGT/0411057.

\bibitem{CO}
P.M. Cohn.
\newblock {\em Free Rings and their Relations, second edition}.
\newblock Number~19 in London Math. Soc Monographs. Academic Press, London,
  1985.

\bibitem{ConT2}
J.~Conant and P.~Teichner.
\newblock Grope cobordism and feynman diagrams.
\newblock {\em Math. Annalen}, 328(1-2):135--171, 2004.

\bibitem{ConT1}
J.~Conant and P.~Teichner.
\newblock Grope cobordism of classical knots.
\newblock {\em Topology}, 43(1):119--156, 2004.

\bibitem{CuSh}
Marc Culler and Peter~B. Shalen.
\newblock Singular surfaces, mod 2 homology, and hyperbolic volume, ii.
\newblock preprint arXiv/math.GT/0701666.

\bibitem{Dw}
W.~Dwyer.
\newblock Homology, massey products and maps between groups.
\newblock {\em Journal of Pure and Applied Algebra}, 6:177--190, 1975.

\bibitem{FQ}
M.H. Freedman and F.~Quinn.
\newblock {\em The Topology of 4-manifolds}.
\newblock Number~39 in Princeton Math. Series. Princeton University Press,
  Princeton, New Jersey, 1990.

\bibitem{Ha2}
S.~Harvey.
\newblock Homology cobordism invariants of 3-manifolds and the
  cochran-orr-teichner filtration of the link concordance group.
\newblock preprint, http://front.math.ucdavis.edu/math.GT/0609378.

\bibitem{HS}
P.J. Hilton and U.~Stammbach.
\newblock {\em A Course in Homological Algebra}, volume~4 of {\em Graduate
  Texts in Mathematics}.
\newblock Springer-Verlag, New York, 1971.

\bibitem{Ho}
J.~Howie.
\newblock Free subgroups in groups of small deficiency.
\newblock {\em J. Group Theory}, 1:95--112, 1998.

\bibitem{HoS}
J.~Howie and H.R. Scheebeli.
\newblock Homological and topological properties of locally indicable groups.
\newblock {\em Manuscripta Math.}, 44 (1--3):71--93, 1983.

\bibitem{Ka2}
Akio Kawauchi.
\newblock Almost identical imitations of {$(3,1)$}-dimensional manifold pairs.
\newblock {\em Osaka J. Math.}, 26(4):743--758, 1989.

\bibitem{Ka1}
Akio Kawauchi.
\newblock An imitation theory of manifolds.
\newblock {\em Osaka J. Math.}, 26(3):447--464, 1989.

\bibitem{Kr2}
V.S. Krushkal.
\newblock On the relative slice problem and and four-dimensional topological
  surgery.
\newblock {\em Math.Annalen}, 315(3):363--396, 1999.

\bibitem{Kr1}
V.S. Krushkal.
\newblock Dwyer's filtration and the topology of four-manifolds.
\newblock {\em Math.Res.Letters}, 10(2-3):247--251, 2003.

\bibitem{KT}
V.S. Krushkal and P.~Teichner.
\newblock Alexander duality, gropes and link homotopy.
\newblock {\em Geom.Top.}, 1:51--69, 1997.

\bibitem{La3}
M.~Lackenby.
\newblock Detecting large groups.
\newblock math.GR/0702571.

\bibitem{La1}
M.~Lackenby.
\newblock Large groups, property (tau) and the homology growth of subgroups.
\newblock math.GR/0509036.

\bibitem{La2}
M.~Lackenby.
\newblock New lower bounds on subgroup growth and homology growth.
\newblock math.GR/0512261.

\bibitem{Lu1}
Alexander Lubotzky and Dan Segal.
\newblock {\em Subgroup growth}, volume 212 of {\em Progress in Mathematics}.
\newblock Birkh\"auser Verlag, Basel, 2003.

\bibitem{FT2}
M.H.Freedman and P.~Teichner.
\newblock 4-manifold topology ii: Dwyer's filtration and surgery kernels.
\newblock {\em Inventiones Math.}, 122:531--557, 1995.

\bibitem{Ru1}
Daniel Ruberman.
\newblock Seifert surfaces of knots in {$S\sp 4$}.
\newblock {\em Pacific J. Math.}, 145(1):97--116, 1990.

\bibitem{SW}
Peter~B. Shalen and Philip Wagreich.
\newblock Growth rates, {$Z\sb p$}-homology, and volumes of hyperbolic
  {$3$}-manifolds.
\newblock {\em Trans. Amer. Math. Soc.}, 331(2):895--917, 1992.

\bibitem{St}
J.~Stallings.
\newblock Homology and central series of groups.
\newblock {\em Journal Algebra}, 2:170--181, 1965.

\bibitem{St2}
J.~Stallings.
\newblock Surfaces in three-manifolds and non-singulat equations in groups.
\newblock {\em Math.Zeit.}, 184:1--17, 1983.

\bibitem{Str}
R.~Strebel.
\newblock Homological methods applied to the derived series of groups.
\newblock {\em Comment. Math.}, 49:302--332, 1974.

\bibitem{CH1}
T.Cochran and S.~Harvey.
\newblock Homology and derived series of groups.
\newblock {\em Geom. Topol.}, 9:2159--2191, 2005.

\bibitem{T}
P.~Teichner.
\newblock Knots, von neumann signatures and grope cobordism.
\newblock {\em Proc. of the ICM, Beijing}, Vol.II: Invited Lectures:437--446,
  2002.

\end{thebibliography}
